\documentclass[letterpaper,11pt]{article}

\title{The Algebraic Complexity of Maximum Likelihood Estimation for Bivariate Missing Data}

\author{Serkan Ho\c{s}ten \\
{\small  Department of Mathematics, San Francisco State University} \\
Seth Sullivant
 \\  {\small Department of Mathematics, Harvard University} }

\date{}

\usepackage{latexsym,array,delarray,amsthm,amssymb}%eufrak
\usepackage[sumlimits]{amsmath}
\usepackage[dvips]{graphicx}
\theoremstyle{plain}
\newtheorem{thm}{Theorem}[section]
\newtheorem{lemma}[thm]{Lemma}

\theoremstyle{definition}

\theoremstyle{remark}

\newcommand{\rr}{\mathbb{R}}
\newcommand{\cc}{\mathbb{C}}

\begin{document}
\maketitle

\begin{abstract}
We study the problem of maximum likelihood estimation for general patterns of bivariate missing data for normal and multinomial random variables, under the assumption that the data is missing at random (MAR).  For normal data, the score equations have nine complex solutions, at least one of which is real and statistically significant.  Our computations suggest that the number of real solutions is related to whether or not the MAR assumption is satisfied.  In the multinomial case, all solutions to the score equations are real and the number of real solutions grows exponentially in the number of states of the underlying random variables, though there is always precisely one statistically significant local maxima.   
\end{abstract}

%%%%%%%%%%%%%%%%%%%%%%%%%%%%%%%%%
%%%%%%%%%%%%%%%%%%%%%%%%%%%%%%%%%
%%%%%%%%%%%%%%%%%%%%%%%%%%%%%%%%%
%%%%%%%%%%%%%%%%%%%%%%%%%%%%%%%%%

\section{Introduction}

A common problem in statistical analysis is dealing with missing covariates in some of the replicates of multivariate data.   A typical instance arises during longitudinal studies in the social and biological sciences, when participants may miss appointments or drop out of the study altogether.  Over very long term studies nearly all replicates will involve some missing data, so it is usually impractical to throw out replicates with missing covariates.  Furthermore, the underlying cause for the censoring (e.g.~a subject dies) might play an important role in inference with the missing data that will lead to false conclusions in the complete case analysis.  Thus, specialized techniques are needed in the setting where some of the data is missing.  A useful reference for this material is \cite{Little2002}, from which we will draw notation and definitions.

In this paper, we undertake an algebraic study of maximum likelihood estimation for general patterns of bivariate missing data, under the assumption that the data is \emph{missing at random} (MAR) \cite{Little2002}.  This implies, in particular, that the censoring mechanism does not affect the maximization of the likelihood function with respect to the underlying parameters of the model, and thus the nonresponse is ignorable.

Let $Y_1, \ldots, Y_n$ be i.i.d. replicates with $d$-covariates $X_1, \ldots, X_d$.  Let $\mathcal{M}$ be a parametric model for the joint distribution of the $X_i$.   
Let $M$ be the $d \times n$ $0/1$-matrix that is the indicator function for the missing entries of the $Y_j$; that is $M_{ij} = 1$ if and only if the $i$th covariate of $Y_j$ is missing.  Under the assumption that the data are missing at random (the missing data mechanism does not depend on the values of the missing data), we can perform likelihood inference with respect to an underlying model for the way that the fully observed data was generated.  In particular, if $f(y | \theta)$ is the probability density function of random variables $X$, we have the log-likelihood function for the observed data:

$$ \ell(\theta | Y, M) =   \sum_{j =1}^n  \log f(Y_j = y_j  | \theta, M),$$
where $f(Y_j = y_j |  \theta, M)$ denotes the marginal probability of observing $Y_j = y_j$ with appropriate entries of $y_j$ censored
$$f(Y_j = y_j | \theta, M)  = \int_{X_i | M_{ij} = 1}   f(X_{\rm obs} = y_{\rm obs}, X_{\rm mis} = x_{\rm mis}| 
\theta )  d x_{\rm mis}. $$
We wish to find the parameter values $\hat{ \theta}$ that maximize this likelihood function.

Our focus in this paper is on the case when $d = 2$.  With a general pattern of missing data in the bivariate case, we assume that our data comes in the following form.  There are $n$ complete cases where we obtain a $2$-dimensional vector $Y_i$. There are $r$ cases where we only obtain variable $X_1$, and $s$ cases where we only obtain variable $X_2$.  We denote these by $Z_i$ and $W_i$, respectively.  The log-likelihood function becomes
$$\ell(\theta; y,w,z)  = \sum_{j =1}^n  \log  f(Y_j = y_j | \theta)  + \sum_{j = 1}^r  \log f(Z_j = z_j | \theta)  + \sum_{j = 1}^s  \log f(W_j = w_j | \theta)$$
and our goal is to maximize this function.

One approach to determining the maximum likelihood estimate uses computational algebraic geometry.  The connections between maximum likelihood estimation and algebraic geometry was first extensively studied in \cite{Catanese2004}.  A basic fact is that, if the critical equations (score equations) are rational functions of the parameters and the data, then the number of complex solutions to the critical equations is constant for generic (i.e. almost all) data.  This fixed number is called the \emph{maximum likelihood degree} (ML-degree for short) of the model.  The ML-degree is an intrinsic complexity measure of the score equations. In this paper, we compute the ML-degree in the bivariate missing data problem for Gaussian random variables and for multinomial random variables.

The outline of this paper is as follows.  In Section \ref{sec:normal} we focus on the case where $(X_1,X_2)$ have a jointly normal distribution.  We show that the ML-degree in this case is nine.  
Our simulations show that if the data is indeed generated from bivariate normal distributions, and the censoring mechanism is MCAR or MAR,   then
there is a unique real solution to the score equations, which is a local maximum.  
On the other hand, we also present examples 
of data, where either the model or the missing data mechanism are misspecified, where there can be two statistically relevant local maxima. The possible existence of multiple maxima is important to take into account when using the EM-algorithm to find the maximum likelihood estimate.  In Section \ref{sec:discrete} we focus on the discrete case, where $(X_1, X_2)$ have a jointly multinomial distribution.  In this setting, we give a combinatorial formula for the ML-degree.

%%%%%%%%%%%%%%%%%%%%%%%%%%%%%%%%%n
%%%%%%%%%%%%%%%%%%%%%%%%%%%%%%%%%
%%%%%%%%%%%%%%%%%%%%%%%%%%%%%%%%%
%%%%%%%%%%%%%%%%%%%%%%%%%%%%%%%%%
\section{Bivariate Normal Random Variables}\label{sec:normal}

We assume that $X = (X_1, X_2)  \sim \mathcal{N}(\mu, \Sigma)$ where $E[X] = \mu = (\mu_1, \mu_2)$
and 
$\Sigma = \left[ \begin{array}{cc} \sigma_{11} & \sigma_{12} \\ \sigma_{12} & \sigma_{22} \end{array} \right]$ is the covariance matrix.  Then we have $Z_j \sim \mathcal{N}(\mu_1, \sigma_{11})$ 
for  $j=1, \ldots, r$ and $W_j \sim \mathcal{N}(\mu_2, \sigma_{22})$ for $j=1, \ldots, s$.  Ignoring constants the log-likelihood function  is equal to
\begin{multline*} 
\ell(\mu, \Gamma |  y,w,z) = -\frac{1}{2}n \log (\det \Sigma) - \frac{1}{2} \big(\sum_{j=1}^n (Y_j - \mu)^t  \Sigma^{-1} (Y_j- \mu) \big) \\
- \frac{1}{2}r \log(\sigma_{11}) - \frac{1}{2\sigma_{11}} \sum_{j=1}^r (Z_j - \mu_1)^2 
 - \frac{1}{2}s \log(\sigma_{22}) - \frac{1}{2\sigma_{22}} \sum_{j=1}^s (W_j - \mu_2)^2
\end{multline*}
It is more convenient to use the entries of 
$\Gamma := \Sigma^{-1} = \left[ \begin{array}{cc} \gamma_{11} & \gamma_{12} \\ \gamma_{12} & \gamma_{22} \end{array} \right]$ in our computations. With this substitution, we get the identities
$$\sigma_{11} = \frac{\gamma_{22}}{\det \Gamma}, \,\,\,   \sigma_{22} = \frac{\gamma_{11}}{\det \Gamma}, 
\,\,\, \mbox{and} \,\,\,  \sigma_{12} = \frac{-\gamma_{12}}{\det \Gamma}. $$
In the computations below we will also use a bar over a quantity to denote its average. 
The log-likelihood function becomes
 \begin{multline}
\frac{1}{2} (n+ r + s)  \log (\det \Gamma) -\frac{1}{2}r \log \gamma_{22} - \frac{1}{2}s \log \gamma_{11} \\
  - \frac{n}{2} \big[(\overline{Y_1^2} - 2 \mu_1 \overline{Y_1} + \mu_1^2) \gamma_{11} 
 + 2 (\overline{Y_1Y_2} - (\overline{Y_1} \mu_2 + \overline{Y_2} \mu_1) + \mu_1 \mu_2) \gamma_{12} 
 + (\overline{Y_2^2} - 2 \mu_2 \overline{Y_2} + \mu_2^2) \gamma_{22}\big] \\
 - \frac{r}{2} \frac{\det \Gamma}{\gamma_{22}} 
 (\overline{Z^2} - 2 \mu_1 \overline{Z} + \mu_1^2)
- \frac{s}{2} \frac{\det \Gamma}{\gamma_{11}} 
 (\overline{W^2} - 2 \mu_2 \overline{W} + \mu_2^2)
\end{multline}
The critical equations for $\ell(\mu, \Gamma; y,z,w)$ are:

\begin{eqnarray}\label{gaussian-critical}
0 = \frac{\partial \ell}{\partial \mu_1}  & =  & n \big[ (\overline{Y_1} - \mu_1)\gamma_{11} + (\overline{Y_2} - \mu_2)\gamma_{12} \big] + r \frac{\det \Gamma}{\gamma_{22}}(\overline{Z} - \mu_1)  \nonumber \\
0 =  \frac{\partial \ell}{\partial \mu_2} & =  & n \big[ (\overline{Y_2} - \mu_2)\gamma_{22} + (\overline{Y_1} - \mu_1)\gamma_{12} \big] + s \frac{\det \Gamma}{\gamma_{11}}(\overline{W} - \mu_2)  \nonumber \\
0 = \frac{\partial \ell}{\partial \gamma_{11}} & = &  \frac{1}{2}(n+r+s) \frac{\gamma_{22}}{\det \Gamma} - \frac{1}{2}\frac{s}{\gamma_{11}} - \frac{n}{2} (\overline{Y_1^2} - 2 \mu_1 \overline{Y_1} + \mu_1^2)  \nonumber \\
 &  &  - \frac{r}{2} (\overline{Z^2} - 2 \mu_1 \overline{Z} + \mu_1^2) - \frac{s}{2} \frac{\gamma_{12}^2}{\gamma_{11}^2}
 (\overline{W^2} - 2 \mu_2 \overline{W} + \mu_2^2) \nonumber \\
 0 =  \frac{\partial \ell}{\partial \gamma_{22}} &  = &  \frac{1}{2}(n+r+s) \frac{\gamma_{11}}{\det \Gamma} - \frac{1}{2}\frac{r}{\gamma_{22}} - \frac{n}{2} (\overline{Y_2^2} - 2 \mu_2 \overline{Y_2} + \mu_2^2) \nonumber \\
 &  &  - \frac{s}{2} (\overline{W^2} - 2 \mu_2 \overline{W} + \mu_2^2) - \frac{r}{2} \frac{\gamma_{12}^2}{\gamma_{22}^2}
 (\overline{Z^2} - 2 \mu_1 \overline{Z} + \mu_1^2) \nonumber \\
0 = \frac{\partial \ell}{\partial \gamma_{12}}  & = &  (n+r+s) \frac{\gamma_{12}}{\det \Gamma} 
 - n  (\overline{Y_1Y_2} - (\overline{Y_1} \mu_2 + \overline{Y_2} \mu_1) + \mu_1 \mu_2) \nonumber \\
  &   &  + r \frac{\gamma_{12}}{\gamma_{22}}
 (\overline{Z^2} - 2 \mu_1 \overline{Z} + \mu_1^2)  + s \frac{\gamma_{12}}{\gamma_{11}}
 (\overline{W^2} - 2 \mu_2 \overline{W} + \mu_2^2)
\end{eqnarray}

%\begin{multline} \label{gaussian-critical}
%0  =   \frac{\partial \ell}{\partial \mu_1} = n \big[ (\bar{Y_1} - \mu_1)\gamma_{11} + (\bar{Y_2} - \mu_2)\gamma_{12} \big] + r \frac{\det \Gamma}{\gamma_{22}}(\bar{Z} - \mu_1)  \\
%0 =  \frac{\partial \ell}{\partial \mu_2} = n \big[ (\bar{Y_2} - \mu_2)\gamma_{22} + (\bar{Y_1} - \mu_1)\gamma_{12} \big] + s \frac{\det \Gamma}{\gamma_{11}}(\bar{W} - \mu_2)  \\
%0 = \frac{\partial \ell}{\partial \gamma_{11}} = (n+r+s) \frac{\gamma_{22}}{\det \Gamma} - \frac{s}{\gamma_{11}} - \frac{n}{2} (\bar{Y_1^2} - 2 \mu_1 \bar{Y_1} + \mu_1^2) \\
%- \frac{r}{2} (\bar{Z^2} - 2 \mu_1 \bar{Z} + \mu_1^2) - \frac{s}{2} \frac{\gamma_{12}^2}{\gamma_{11}^2}
% (\bar{W^2} - 2 \mu_2 \bar{W} + \mu_2^2)\\
%0 =  \frac{\partial \ell}{\partial \gamma_{22}} = (n+r+s) \frac{\gamma_{11}}{\det \Gamma} - \frac{r}{\gamma_{22}} - \frac{n}{2} (\bar{Y_2^2} - 2 \mu_2 \bar{Y_2} + \mu_2^2) \\
%- \frac{s}{2} (\bar{W^2} - 2 \mu_2 \bar{W} + \mu_2^2) - \frac{r}{2} \frac{\gamma_{12}^2}{\gamma_{22}^2}
% (\bar{Z^2} - 2 \mu_1 \bar{Z} + \mu_1^2) \\
% 0 = \frac{\partial \ell}{\partial \gamma_{12}} = -2(n+r+s) \frac{\gamma_{12}^2}{\det \Gamma} 
% - n  (\bar{Y_1Y_2} - (\bar{Y_1} \mu_2 + \bar{Y_2} \mu_1) + \mu_1 \mu_2) \\
% + r \frac{\gamma_{12}}{\gamma_{22}}
% (\bar{Z^2} - 2 \mu_1 \bar{Z} + \mu_1^2)  + s \frac{\gamma_{12}}{\gamma_{11}}
 %(\bar{W^2} - 2 \mu_2 \bar{W} + \mu_2^2)
%\end{multline}

\begin{thm} \label{thm:gaussian}
The ML-degree of the bivariate normal missing data problem is equal to  nine, and
at least one of the critical solutions  to (\ref{gaussian-critical}) is real. Moreover, at least
one such real critical solution is a local maximum in the statistically relevant parameter space.
\end{thm}

\begin{proof} The theorem follows from a general principle about the number of complex 
solutions to  a system of polynomial equations with parametric coefficients. Namely, if such a
system has $N < \infty$ complex solutions (counted with multiplicity)  for a "random" choice of parameter values then other random choices of parameter values will also produce $N$ complex
solutions. Here we sketch a proof of this statement. Suppose $I$ is an ideal in
$\cc(p_1, \ldots, p_k)[x_1, \ldots, x_t]$, the ring of polynomials in the indeterminates $x_1, \ldots, x_n$
with coefficients from the field of rational functions in $p_1, \ldots, p_k$ over $\cc$. Pick any term order
and compute a Gr\"obner basis $G$ of $I$ with respect to this term order. Now let $U$ be the Zariski
open set  in $\cc^k$ such that no denominator of the coefficients and no initial coefficient of the polynomials 
encountered during the Buchberger algorithm that produces $G$ vanish on any point in $U$. If
${\bar p} \in U$ then both the initial ideal of $I$ and that of $I({\bar p})$ will have the same set of
standard monomials: these are the monomials that no initial term in $G$ and $G({\bar p})$, respectively, divide. It is a well-known result that 
$I({\bar p})$ has $N < \infty$ complex solutions (counted with multiplicity) if and only if 
the number of such standard monomials is $N$.   This implies that for all ${\bar q} \in U$ the ideal
$I({\bar q})$ will have $N$ complex solutions.

Now in the setting of the critical equations (\ref{gaussian-critical}) let  $J$ be the ideal generated
by the five polynomials obtained by clearing the denominators in (\ref{gaussian-critical}). Furthermore,
let $K$ be the ideal generated by the product of these cleared denominators. Then the ML-degree
we are after is the number of complex solution of $I = J  : K$. A random choice of $n, r, s$ and data 
vectors $y_1, \ldots, y_n$, $z_1, \ldots, z_r$, and $w_1, \ldots, w_s$, and a quick computation
in Singular shows that $I(n,r,s,y,w,z)$ has nine complex solutions. Our discussion above
implies that the ML-degree of the bivariate normal missing data problem is nine. Since complex
solutions to real polynomial equations come in complex conjugate pairs, at least one must be a real solution. 
  
We can also see directly that there must be at least one real local maximum inside the statistically relevant parameter space $\rr^2 \times PD_2$ (where $PD_2$ denotes the space of $ 2 \times 2$ positive definite matrices).  To see this, note that if any parameter has a large absolute value the log-likelihood function tends to $-\infty$. Similarly, if the $\Sigma$ parameters approach the boundary of the positive definite cone the log-likelihood function tends to $- \infty$.  Thus, the log-likehood function must have a local maximum inside $\rr^2 \times PD_2$.  
\end{proof}

How many of the nine complex solutions in Theorem \ref{thm:gaussian}  can be real? We know that
at least one is, but is it possible that there are three, five, seven, or nine?  For various choices of the data parameters, we have observed that all of these values are possible.  A more surprising fact is that the number of real solutions seems to be indicative of how well-specified the MAR assumption is.  Here is a summary of the observations that emerge from our computations 
for which we have use Mathematica, Maple, and Singular \cite{GPS05}.  We describe the separate cases in more detail in the paragraphs following the list.

\begin{enumerate}
\item  When the data was generated from a Gaussian or uniform distribution  and the missing data mechanism was MCAR (missing completely at random) or MAR, we consistently observe exactly one real critical point, which is necessarily a local maximum.
\item  When the data was generated from a Gaussian distribution and the missing data mechanism was NMAR (not missing at random), we consistently observed three real critical points, all of which are in the statistically relevant region ($\rr^2 \times PD_2$), of which two were local maxima.
\item  When the joint distribution of $Y$ and the marginal distributions of $W$ and $Z$ were unrelated to each other by a natural censoring mechanism, we observed seven real critical points, of which three were in the statistically relevant region, and two were statistically relevant local maxima.
\item  When the twelve sufficient statistics ($n,r,s, \overline{Y_1}, \ldots$) were generated randomly (without regard to an underlying distribution) we observed nine real critical points. 
\end{enumerate}

Of course, we could not test all possible scenarios for the above data types, and there will always be the possibility that data generated by one of the strategies will have a different number of real solutions than we observed.

When the censoring mechanism was MCAR, we did the censoring in the obvious way, by first generating data from a randomly chosen Gaussian distribution, and then censoring cell entries with the fixed probability $\frac{1}{5}$.  For a more general MAR scenario, we generated data by taking a mixture of the MCAR scenario, with the censoring mechanism that covariate $X_2$ is not observed whenever $X_1 < -1$.  Out of 1000 runs of the MAR scenario 985 cases produced a single 
real solution which is also a statistically relevant maximum.  In fact, both of the above  scenarios consistently had one real solution.

For the NMAR censoring mechanism, we generated data from a random, strongly negatively-correlated Gaussian distribution, and censored covariate $X_i$ when $X_i < -1$.  Out of 1000 sample runs under
this  scenario 765 generated three real solutions, all statistically relevant, with two being local maxima.

For a family of ``wild'' examples, we choose $Y$ and $Z$ to be generated from the same Gaussian distributions with mean $(0,0)$ but $W$ to be generated from   a uniform distribution on the interval $[5,6]$.  We tested this scenario with 1000 sample runs as well, and we observed 831 of them
having seven real solutions, three of  them statistically significant, with two local maxima.  

For the case of randomly generated data without regard to an underlying distribution we also did
1000 sample runs where we observed 134 cases with nine real critical sollutions.

In summary, our computations suggest that the number of real solutions of the critical equations can be a gauge of how well the MAR  assumption fits the data.  For missing data sets with three or more covariates where direct computation of all critical points will not be possible, if the EM-algorithm produces more than one local maximum, this might suggest that one should pay more careful attention to whether or not the MAR assumption makes sense for the data.

%The above example is a "wild" example in the sense that the data is generated in an arbitrary
%manner without any structure in it. It is highly unlikely that a real world missing data
%problem will be such a wild example. For instance, if we actually generate a missing data problem %from $(X_1, X_2)  \sim \mathcal{N}(\mu, \Sigma)$ and find the critical solutions to the score equations
%we hope to get few number of real solutions. In order to test this claim, we have generated
%6,000 instances where $ 50 \leq n \leq 100$ and $5 \leq r,s \leq 30$ were integers chosen uniformly
%at random and where the data vectors, $y,z$ and $w$ were drawn from $\mathcal{N}(\mu, \Sigma)$
%for random choices of $\mu$ and $\Sigma$. In {\em all} of these test runs the score equations 
%had a unique real solution. It is also interesting to observe that if one replaces 
%$\mathcal{N}(\mu, \Sigma)$ with a continuous uniform distribution in these test runs the result
%does not change:  the score equations have a unique real solution.  However, we do see
%a deviation from this behavior if the missing data is not MAR.

%\begin{ex} Please insert the not MAR example.
%\end{ex} 

%%%%%%%%%%%%%%%%%%%%%%%%%%%%%%%%%
%%%%%%%%%%%%%%%%%%%%%%%%%%%%%%%%%
%%%%%%%%%%%%%%%%%%%%%%%%%%%%%%%%%
%%%%%%%%%%%%%%%%%%%%%%%%%%%%%%%%%

\section{Bivariate Discrete Random Variables}\label{sec:discrete}

In this section, we focus on the case where $X_1$ and $X_2$ are discrete multinomial random variables.  We suppose that $X_1 \in \{1,2,\ldots, m\}$ and $X_2 \in \{1,2, \ldots, n\}$.  We give a combinatorial formula the ML-degree which shows that it grows exponentially as a function of $m$ and $n$.  

In the bivariate multinomial case, the data can be summarized by a table of counts $T$ which records the complete cases, and two vectors $R$ and $S$ which record the observations of only $X_1$ and only $X_2$, respectively.  In this multinomial case, we want to estimate the raw probabilities $p_{ij} = P(X_1 = i, X_2 = j)$.  The log-likelihood function becomes
\begin{equation}\label{eqn:bi}   \ell(R,S,T; p) = \sum_{i =1}^m \sum_{j =1}^n   t_{ij} \log p_{ij}  + \sum_{i =1}^m  r_i \log p_{i+}  + \sum_{j = 1}^n s_j \log p_{+j}.\end{equation}
We want to find $p$ that maximizes $\ell(R,S,T; p)$ subject to $p \geq 0$ and $p_{++} = 1$.

\begin{thm}\label{thm:basicdiscrete}
The ML-degree of the bivariate multinomial missing data problem is equal to the number of bounded regions in the arrangement of hyperplanes $\{p_{ij} = 0, p_{i+} = 0, p_{+j} = 0 \, | \,  i \in [m], j \in [n] \}$ inside the hyperplane $p_{++} = 1$.   Every solution to the score equations for (\ref{eqn:bi}) is real.  
For generic $R,S,T$ there is exactly one nonnegative critical point, and it is a local maximum.
\end{thm}

\begin{proof}
Maximizing the product of linear forms has a standard formula for the ML-degree as the number of bounded regions in the arrangement defined by these linear  forms \cite{Catanese2004}.  
Each bounded region contains precisely one critical solution which is real.
Furthermore, since all the coordinate probability functions are linear in the parameters, the objective function is convex  so there is exactly one nonnegative criticial point that must be a local maximum.
\end{proof}

From Theorem \ref{thm:basicdiscrete} we see that to calculate the ML-degree we need to count the number of bounded regions in a hyperplane arrangement.  The remainder of this section is devoted to performing this count.  First we provide some definitions which allow us to state  Theorem \ref{thm:count}.  Then we proceed with the proof  in a number of steps.

For integers $k$ and $l$, the \emph{Stirling numbers of the second kind} are the numbers
$$S(l,k)  =  \frac{1}{k!}  \sum_{i = 0}^k (-1)^{k -i}  {  k  \choose i}  i^l.$$
The \emph{negative index poly-Bernoulli numbers} are the numbers:
$$B(l,k)  =   \sum_{i =0}^l  (-1)^{l - i}  i!  S(l,i) (i + 1)^k. $$

\begin{thm}\label{thm:count}
The ML-degree of the bivariate multinomial $m \times n$ missing data problem is
\begin{equation}\label{eqn:single}
ML(m,n) =  \sum_{k = 0}^m  \sum_{l = 0}^n  (-1)^{m+n - k - l}  { m \choose k}  {n \choose l} B(m - k, n - l).\end{equation}
\end{thm}

For small values of $m$, we can explicitly work out formulas for this ML-degree.  In particular, one can show that $ML(2,n) =  2^{n+1} - 3$.  Since the ML-degree is monotone as a function of $m$ and $n$, this shows that the ML-degree in the bivariate discrete case is exponential in the size of the problem.
 Let
$$S = \{p_{ij}, \, | \, i \in [m] \cup \{+\}, j \in [n] \cup \{+\}  \} \setminus \{p_{++} \}$$ 
be the set of all hyperplanes in the hyperplane arrangement that determines the ML-degree.
Specifying a (possibly empty) region of the arrangement amounts to choosing a partition $S = N \cup P$.      The resulting open region on the hyperplane $p_{++} = 1$ consists of all matrices $p$ such that $p_{ij} < 0$ if $p_{ij} \in N$ and $p_{ij} > 0$ if $p_{ij} \in P$ and $\sum_{i,j} p_{ij} = 1$.   We denote this set of matrices by $\mathcal{M}(N,P)$.  Our goal is characterize and count the partitions $N \cup P$ such that $\mathcal{M}(N,P)$ is nonempty and bounded.  We prove a sequence of results classifying the type of subconfigurations that can appear in $N$ and $P$.

\begin{lemma}\label{lem:minor}
Let $i, k \in [m]$ with $i \neq k$ and $j,l \in [n]$ with $j \neq l$.  Suppose that $p_{ij}, p_{kl} \in N$ and $p_{il}, p_{kj} \in P$.  Then if $\mathcal{M}(N,P)$ is nonempty it is unbounded.
\end{lemma}

\begin{proof}
Let $e_{ij}$ denote the $m \times n$ matrix with a one in the $ij$ position and zeros elsewhere.  Suppose that $p \in \mathcal{M}(N,P)$.  Then $p + a ( e_{il} + e_{kj} - e_{ij} - e_{kl}) \in \mathcal{M}(N,P)$ for all $a > 0$ since adding $a ( e_{il} + e_{kj} - e_{ij} - e_{kl})$ does not change the sign of any entry of $p$ nor does it change any of the margins $p_{i+}$ of $p_{+j}$.    Thus $\mathcal{M}(N,P)$ contains matrices with arbitrarily large entries and is unbounded.
\end{proof}

Let $N'  = N \cap \{p_{ij} \, | \, i \in [m], j \in [n]  \}$ and $P'  = P \cap \{p_{ij} \, | \, i \in [m], j \in [n]  \}.$  A partition $\lambda = (\lambda_1, \ldots, \lambda_m) $ is a nonincreasing sequence of nonnegative integers.  The length of $\lambda$ is $m$ (we allow zeros in the partition).

\begin{lemma}\label{lem:tableau}
Suppose that $\mathcal{M}(N,P)$ is nonempty and bounded.  There exists a permutation $\sigma$ of the rows and columns of $p$ and a partition $\lambda$ such that 
$$\sigma(N') = \{p_{ij}  \, | \, j \leq  \lambda_i \}.$$
The same is true for $P'$ and for every rectangular submatrix of $p$.
\end{lemma}

\begin{proof}
After permuting rows we may assume that the number of elements in row $i$, $\lambda_i$, is a nonincreasing sequence.  Permuting the columns we may suppose that the only elements of $N'$ in the first row of $p$ are $p_{11}, \ldots, p_{1\lambda_1}$.  Permuting columns further, we may assume that the elements in the second row are of the form $p_{21}, \ldots, p_{2 \lambda_2}$ with $\lambda_2 \leq \lambda_1$.  There could not be any element of the form $p_{2j} \in N'$ with $j > \lambda_1$ because otherwise there would be more entries in row two than row one or $N'$ would contain $p_{1\lambda_1}, p_{2j}$ and $P'$ would contain $p_{1j}, p_{2\lambda_1}$ which violates Lemma \ref{lem:minor}.  Repeating the argument for each row shows that $\mathcal{M}(N,P)$ can be put into partition form.
\end{proof}

\begin{lemma}\label{lem:P}
Suppose that $\mathcal{M}(N,P)$ is nonempty and bounded.  Then $p_{i+}, p_{+j} \in P$ for all $i$ and $j$.
\end{lemma}

\begin{proof}
Suppose that $\mathcal{M}(N,P)$ is nonempty and $N$ contains, say, $p_{+1}$.  We will show $\mathcal{M}(N,P)$ is unbounded.  To do this, it suffices to show that there exist points on the boundary of $\mathcal{M}(N,P)$ with coordinates of arbitrarily large absolute values.  Furthermore, we will assume that $\mathcal{M}(N,P)$ is bounded (so that we can make liberal use of Lemmas \ref{lem:tableau} and \ref{lem:minor}) and derive a contradiction.  The boundary of $\mathcal{M}(N,P)$ is described by allowing the strict inequalities to become weak inequalties.  There are four cases to consider.

 {\em Case 1}.  Suppose that there is no $i$ such that $p_{i+} \in N$.  After permuting columns and rows we may suppose that $p_{+j} \in N$ if and only if $j \in [k]$. If $\mathcal{M}(N,P)$ is to be nonempty, we must have $k < m$.

After permuting row and columns in such a way that the set of the first $k$ columns is mapped to itself, we may suppose that the set of variables in $N$ belonging to the submatrix $p[1,m; 1, k]$ is in partition form, according to Lemma \ref{lem:tableau}.  If $\mathcal{M}(N,P)$ is to be nonempty, it must be the case that $p_{1j} \in N$ for all $j \in [k]$ since the first row is the longest row of the tableau.  As $p_{i+} \in P$, there must exist $p_{1l} \in P$ with $l > k$.  Then consider the matrix $p'$ with $p'_{11} = -a$, $p_{1j} = a + 1$ and $p_{ij} = 0$ for all other $i,j$.  This matrix satisfies all requirements to belong to the boundary of $\mathcal{M}(N,P)$.  Letting $a$ tend to infinity shows that $\mathcal{M}(N,P)$ is unbounded, a contradiction.
 
\smallskip
For the remaining three cases, we assume that there exists some $i$ and $j$ such that $p_{i+}, p_{+j} \in N$.  After permuting rows and columns we may suppose there is $k < m$ and $l < n$ such that $p_{i+} \in N$ if and only if $i \in [k]$ and $p_{+j} \in N$ if and only if $j \in [l]$.
 
{\em Case 2}.
Suppose that there is a $p_{ij} \in N$ with $i\in [k]$ and $j \in [l]$ and a $p_{i'j'} \in P$ with $i' \in [k+1, m]$ and $j' \in [l+1, n]$.  Then the matrix $p'$ with $p_{ij} = -a$, $p_{i'j'} = a+1$ and all other entries equals satisfies the requirements to belong to the boundary of $\mathcal{M}(N,P)$.  Letting $a$ tend to infinity shows that $\mathcal{M}(N,P)$ is unbounded, a contradiction.

 {\em Case 3}.
Suppose that $p_{ij} \in P$ for all $i \in [k]$ and $j \in [l]$.  Since $\mathcal{M}(N,P)$ is nonempty, and  $p_{i+} \in N$ for all $i \in [k]$, we can find, for each $i \in [k]$, a $j \in [l+1, n]$ such that $p_{ij} \in N$.  As $\mathcal{M}(N,P)$ is bounded, this implies that we can permute rows and columns of the matrix $p$, so that $p[1,k;l+1,n]$ is mapped into itself and so that this submatrix, intersected with $N$ is of tableau form.  With these assumptions, we must have $p_{il+1} \in N$ for all $i \in [k]$.  Since $p_{+,l+1} \in P$, there must exist $p_{i'l+1} \in P$ with $i' \in [k+1,m]$.   Now consider the matrix $p'$ with $p'_{1l+1} = -a$, $p'_{
i'l+1} = a+1$ and all other entries equal to zero.  This matrix satisfies all requirements for belonging to the boundary of $\mathcal{M}(N,P)$ but as $a$ tends to infinity shows that $\mathcal{M}(N,P)$ is unbounded.

{\em Case 4}.  
Suppose that $p_{ij} \in N$ for all $i \in [k+1, m]$ and $j \in [l+1, n]$.  This is equivalent to saying that for all $p_{ij} \in P$, $p_{i+}$ and $p_{+j}$ are not simultaneously in $P$.  If we permute rows and columns of $p$ so that $P$ is in tableau form, this condition is equivalent to saying that there is a $p_{i'j'} \in P$ such that $p_{i'+1j'+1} \notin P$ and none of  $p_{i+}$ nor $p_{+j}$ are in $P$ for $i \leq i'$ and $j \leq j'$.  (Note that one of $i'$ or $j'$ might be zero, which will work fine in the following argument.)
Then for any matrix $p \in \mathcal{M}(N,P)$ we have
\begin{eqnarray*}
 0  & >  &   \sum_{i = 1}^{i'}  p_{i+}  +  \sum_{j = 1}^{j'} p_{+j} \\
  &  =  &     2 \sum_{i = 1}^{i'}  \sum_{j = 1}^{j'} p_{ij}  + 
                 \sum_{i = i'+1}^m  \sum_{j = 1}^{j'} p_{ij} +
                 \sum_{i = 1}^{i'}  \sum_{j = j' +1}^n p_{ij}
\end{eqnarray*}
The expression at the end of this equation involves the sum, with positive coefficients, of all $p_{ij} \in P$.  Since the $p_{ij}$ in the sum with $p_{ij} \in N$ all occur with coefficient $1$, and since $p_{++} = 1$, we deduce that this sum must be strictly greater than $1$.  Thus $\mathcal{M}(N,P)$ must be empty.
\end{proof}

\begin{lemma}\label{lem:bounded}
Let $\lambda$ be a partition of length $m$ such that $\lambda_i \leq n-1$ for all $i$, and $\lambda_m = 0$.  Let $N(\lambda) = \{ p_{ij} \, | \, j \leq \lambda_i \}$ and $P(\lambda) = S \setminus N(\lambda)$.  Then $\mathcal{M}(N(\lambda),P(\lambda))$ is nonempty and  bounded.
\end{lemma}

\begin{proof}
To show that $\mathcal{M}(N(\lambda),P(\lambda))$ is nonempty amounts to showing that there is a table $p$ with nonzero entries that satisfies all the constraints $p_{ij} < 0$ if $p_{ij} \in N(\lambda)$, $p_{ij} > 0 $ if $p_{ij} \in P(\lambda) $ and $p_{++} = 1$.  To this end, let $\epsilon >0$ be a small real number.  Define the matrix $p(\epsilon)$ by the following rules:
$$p(\epsilon)_{ij}  = \left\{ \begin{array}{cl}
-\epsilon & \mbox{if } p_{ij} \in N(\lambda)  \\
\epsilon &  \mbox{if }  p_{ij} \in P(\lambda) \mbox{ and }  i < m, j < n  \\
m \epsilon & \mbox{if }  i = m, j < n \\
n \epsilon & \mbox{if }  i < m, j = n  \\
1 - (3mn -2m - 2n + 1 - 2 \sum_k \lambda_k) \epsilon  & \mbox{if } i = m, j = n
\end{array} \right.$$ 
By construction, $p( \epsilon) \in \mathcal{M}(N,P)$.

Now we show that $\mathcal{M}(N(\lambda), P(\lambda) )$ is bounded.  For each $k \in [m-1]$ with $\lambda_k > 0$ we have
\begin{eqnarray*}
0 & \leq  &  \sum_{i = 1}^k  p_{i+}  +  \sum_{j = 1}^{\lambda_k} p_{+j} \\
 &  =  & 2 \sum_{i = 1}^k  \sum_{j = 1}^{\lambda_k} p_{ij}  + 
                 \sum_{i = k+1}^m  \sum_{j = 1}^{\lambda_k} p_{ij} +
                 \sum_{i = 1}^k  \sum_{j = \lambda_k + 1}^n p_{ij}
\end{eqnarray*}
which implies that 
\begin{eqnarray*}
-\left( \sum_{i = 1}^k  \sum_{j = 1}^{\lambda_k} p_{ij} \right)
 & \leq &
   \sum_{i = 1}^k  \sum_{j = 1}^{\lambda_k} p_{ij}  + 
                 \sum_{i = k+1}^m  \sum_{j = 1}^{\lambda_k} p_{ij} +
                 \sum_{i = 1}^k  \sum_{j = \lambda_k + 1}^n p_{ij} \\
   &  \leq & \sum_{i = 1}^m \sum_{j = 1}^n p_{ij}  = 1.
\end{eqnarray*}
Since $p_{ij} \in N(\lambda)$ whenever $i \in [k]$ and $j \in [\lambda_k]$, we deduce that 
$$ -1 \leq \sum_{i = 1}^k  \sum_{j = 1}^{\lambda_k} p_{ij} \leq 0 $$
and thus $-1 \leq p_{ij} \leq 0$.  Since every $p_{ij} \in N(\lambda)$ belongs to such a sum for some $k$, we see that $p_{ij}$ is bounded for all $p_{ij} \in N(\lambda)$.  This implies that $p_{ij}$ is bounded for all $p_{ij} \in P(\lambda)$ as well, since, $p_{++} = 1$.  Thus, $\mathcal{M}(N(\lambda), P(\lambda))$ is bounded. 
\end{proof}

To finish the proof, we use a result from the Master's thesis of Chad Brewbaker \cite{Brewbaker2005}, that counts a family of $0/1$ matrices that are closely related to the set $N,P$ that have $\mathcal{M}(N,P)$ bounded.

\begin{thm}
The number of $0/1$ $m \times n$ matrices $A$ such that no $ 2 \times 2$ submatrix of $A$ is either $\begin{pmatrix}1 & 0 \\ 0 & 1 \end{pmatrix}$ or $\begin{pmatrix}0 & 1 \\ 1 & 0 \end{pmatrix}$ is the negative index poly-Bernoulli number $B(m,n)$.
\end{thm}

The $0/1$ matrices in the theorem are known as \emph{lonesum matrices} because they are the $0/1$ matrices that are uniquely specified by their row and column sums.

\smallskip

\noindent  \emph{Proof of Theorem \ref{thm:count}}
According to Lemmas \ref{lem:minor}, \ref{lem:P}, and \ref{lem:bounded}, we must count sets $N \subset \{p_{ij} \,\, | \, \, i \in [m], j \in [n]  \}$ with certain properties.  Interpreting $N$ as a $0/1$ matrix where $M$ where $M_{ij} = 1$ if $p_{ij} \in N$, we see that we must count the matrices $M$ that do not have any $2 \times 2$ submatrices equal to $\begin{pmatrix}1 & 0 \\ 0 & 1 \end{pmatrix}$ or $\begin{pmatrix}0 & 1 \\ 1 & 0 \end{pmatrix}$.   Furthermore, the fact that no $p_{i+}$ or $p_{+j}$ belongs to $N$ implies that no row or column of $M$ could be all ones (otherwise, we would have, for example, $p_{ij} < 0 $ for all $j$ but $p_{i+} > 0$ which implies that $\mathcal{M}(N,P)$ is empty) .  Because of the fact that each such set $N$ can be rearranged into a partition, and after switching the zeros and ones, this is the same as the number of lonesum $0/1$ $m \times n$ matrices which have all row and column sums positive.  Thus, the number $M(m,n)$ can be obtained from the negative index poly-Bernoulli  numbers $B(m,n)$ by inclusion-exclusion which yields the desired formula (\ref{eqn:single}).
\qed


\begin{thebibliography}{99}

\bibitem{Brewbaker2005}  C.~Brewbaker.  \emph{Lonesum $(0,1)$-matrices and poly-Bernoulli numbers.}  Master's Thesis,  Iowa State University, 2005.

\bibitem{Catanese2004}  F.~Catanese, S.~Ho\c{s}ten, A.~Khetan, and B.~Sturmfels.  The maximum likelihood degree.  \emph{Amer. J. Math.} {\bf  128} (2006), no. 3, 671--697.

\bibitem{GPS05} G.~-M.~Greuel, G.~Pfister, and H.~Sch\"onemann. {\sc Singular} 3.0. A Computer Algebra System for Polynomial Computations. Centre for Computer Algebra, University of Kaiserslautern
(2005). {\tt http://www.singular.uni-kl.de}.

\bibitem{Little2002}  R.~Little and D.~Rubin.  \emph{Statistical Analysis with Missing Data}.  Series in Probability and Statistics, Wiley Interscience, Hoboken, New Jersey, 2002.


\end{thebibliography}
\end{document}